\RequirePackage{fix-cm}
\documentclass[smallextended]{svjour3}       % onecolumn (second format)
\smartqed  % flush right qed marks, e.g. at end of proof
\usepackage{graphicx}
\journalname{Optimization Letters}

\usepackage{algorithm}
\usepackage{algorithmic}
\usepackage{amsmath}

\renewcommand{\span}{\mbox{span}}

\newcommand{\comment}[1]{} % for invisible comments

\DeclareMathOperator{\rank}{rank}

\DeclareMathOperator*{\argmin}{argmin}
\DeclareMathOperator{\sgn}{sgn}

\renewcommand{\dfrac}{\displaystyle \frac} % for better display of fractions

\begin{document}

\title{A Block Lanczos with Warm Start Technique for Accelerating Nuclear Norm Minimization Algorithms%\thanks{Grants or other notes
}


\author{Zhouchen Lin \and
        Siming Wei %etc.
}


\institute{  Zhouchen Lin (corresponding author) \at
              Microsoft Research Asia, 5th Floor, Sigma Building, Zhichun Road \#49, Haidian District, Beijing 100190, P.R. China\\
             \email{zhoulin@microsoft.com}
           \and
              Siming Wei \at
              Zhejiang University\\
%              Tel.: +123-45-678910\\
%              Fax: +123-45-678910\\
              \email{tobiawsm@gmail.com}           %  \\
}

\date{Received: date / Accepted: date}

\maketitle

\begin{abstract}
Recent years have witnessed the popularity of using rank
minimization as a regularizer for various signal processing and
machine learning problems. As rank minimization problems are
often converted to nuclear norm minimization (NNM) problems, they have
to be solved iteratively and each iteration requires computing a
singular value decomposition (SVD). Therefore, their solution
suffers from the high computation cost of multiple SVDs. To
relieve this issue, we propose using the block Lanczos method to
compute the partial SVDs, where the principal singular subspaces obtained
in the previous iteration are used to start the block Lanczos
procedure. To avoid the expensive reorthogonalization in the
Lanczos procedure, the block Lanczos procedure is performed
for only a few steps. Our block Lanczos with warm start (BLWS)
technique can be adopted by different algorithms that solve
NNM problems. We present numerical results
on applying BLWS to Robust PCA and Matrix Completion problems.
Experimental results show that our BLWS technique usually
accelerates its host algorithm by at least two to three times.

\end{abstract}

\section{Introduction}
\label{intro} In recent years, there is a surge of applying rank
minimization as a regularizer to various machine learning and
signal processing problems
\cite{Wright09,Candes08,Wu10,Zhang10-TILT,Zhu10,Min10,Liu10,Peng10,Candes10,Candes09,Meng10,Liu09,Ji10,Batmanghelich10,Karbasi10}.
In the mathematical models of these problems, the rank of some matrix
is often required to be minimized. Typical models are Robust PCA
(RPCA) \cite{Wright09}:
\begin{equation}
\mbox{\bf(RPCA)}\quad \min_{A,E} \rank(A)+\lambda\|E\|_{l_0},\quad
s.t.\quad D=A+E,
\end{equation}
and Matrix Completion (MC) \cite{Candes08}:
\begin{equation}
\mbox{\bf(MC)}\quad \min_{A} \rank(A),\quad s.t.\quad
D=\pi_{\Omega}(A),
\end{equation}
where $\|E\|_{l_0}$ is the number of nonzeros in $E$, $\Omega$ is
the set of indices of known entries in $A$ and $\pi_{\Omega}$ is the
restriction onto $\Omega$. There are variations of RPCA
\cite{Zhou10,Ganesh10} and MC \cite{Candes09}, and there is also a
combination of RPCA and MC \cite{Candes10}.

Due to the effectiveness of rank minimization, many researchers
have proposed various algorithms to solve rank minimization
problems
\cite{Toh09,Cai08,Ma09,Lin09,Ganesh09,Yuan09,Tao09,Candes08}. As
rank minimization problems are usually NP hard, most of them
aim at solving companion convex programs instead, by replacing the
rank function with the nuclear norm $\|\cdot\|_*$, i.e., the sum
of the singular values, and the $l_0$ norm with the $l_1$ norm,
i.e., the sum of the absolute values of the entries. This is
suggested by the fact that the nuclear norm and $l_1$ norm are the
convex envelopes of the rank function \cite{Recht10} and the $l_0$
norm, respectively. Some researchers have proven that for RPCA and
MC problems solving the companion convex program can produce the
same solution to the original problem at an overwhelming
probability \cite{Recht10,Candes09,Candes10}. As a result, solving
a rank minimization problem is often converted into solving a
nuclear norm minimization (NNM) problem, in order to exploit the
efficiency of convex programs.

Whichever of the existing methods that solve the NNM problems is used, one always has to solve the
following subproblem:
\begin{equation}\label{eqn:sub_problem}
A_{i+1}=\argmin_{A}
\varepsilon_i\|A\|_*+\dfrac{1}{2}\|A-W_i\|_{F}^2,
\end{equation}
where $\varepsilon_i$ and $W_i$ change along iteration and $\|\cdot\|_F$ is the Frobenius norm. Cai et al.
\cite{Cai08} proved that the solution to (\ref{eqn:sub_problem})
can be obtained by singular value thresholding:
\begin{equation}\label{eqn:SVT}
A_{i+1}=\mathcal{T}_{\varepsilon_i}(W_i)\equiv
U_i\Theta_{\varepsilon_i}(S_i)V_i^T,
\end{equation}
where $\Theta_\varepsilon(x)=\sgn(x)\max(|x|-\varepsilon,0)$
%\begin{equation}
%\Theta_\varepsilon(x) \doteq \left\{ \begin{array}{ll} x -
%\varepsilon, & \mbox{if }x > \varepsilon, \\ x + \varepsilon, &
%\mbox{if }x < - \varepsilon,
%\\ 0, & \text{otherwise},
%\end{array} \right.
%\end{equation}
is a shrinkage operator and $U_iS_iV_i^T$ is the singular value
decomposition (SVD) of $W_i$. Therefore, it is easy to see that
NNM problems are usually computationally
costly as they require solving SVDs multiple times and an SVD
typically requires $O(p^3)$ operations, where $p=\min(m,n)$ and
$m\times n$ is the size of the matrix. Fortunately, it is apparent
that all the singular values/vectors need not be computed because
the singular values smaller than the threshold $\varepsilon_i$
will be shrunk to zeros hence their associated singular vectors
will not contribute to $A_{i+1}$. This leads to a common practice
in solving NNM problems, namely using
PROPACK \cite{Propack} to compute the partial SVD of $W_i$, where
only those leading singular values that are greater than
$\varepsilon_i$, and their associated singular vectors, are
computed. This significantly brings down the computation
complexity from $O(p^3)$ to $O(rp^2)$, where $r$ is the number of
leading singular values/vectors computed.

Although computing the partial SVD instead already saves the
computation significantly, the $O(rp^2)$ complexity is still too
high for large scale problems. Therefore, any further savings in
the computation are valuable when the problem scale becomes large.
In this paper, we aim at exploiting the relationship between
successive iterations to further bring down the computation cost.
Our technique is called the block Lanczos with warm start (BLWS),
which uses the block Lanczos method to solve the partial SVD and the
block Lanczos procedure is initialized by the principal
singular subspaces of the previous iteration. The number of
steps in the block Lanczos procedure is also kept small. Our BLWS
technique can work in different algorithms for NNM problems. Our numerical tests show that BLWS can speed up its
host algorithm by at least two to three times.

To proceed, we first introduce how the partial SVD is
computed in PROPACK.

\section{The Lanczos Method for the Partial SVD}
\label{sec:Lanczos} PROPACK uses the Lanczos method to compute the partial SVD.
As the method is based on the Lanczos method for partial eigenvalue decomposition (EVD),
we have to start with the partial EVD computation.

The Lanczos method for partial EVD is to find the optimal
leading\footnote{Actually it can also find the tailing
eigen-subspace of $W$.} eigen-subspace of a symmetric matrix $W$
in a Krylov subspace \cite{Golub96}:
\begin{equation}\label{eq:Lanczos}
K(W,q_1,k)=\span\{q_1,Wq_1,\cdots,W^{k-1}q_1\}.
\end{equation}
The orthonormal basis $Q_k$ of $K(W,q_1,k)$ can be efficiently
computed via a so-called Lanczos procedure shown in
Algorithm~\ref{Lanczos}. Accordingly, $W$ can be approximated as
$W\approx Q_k T_kQ_k^T$, where $T_k$ is a tri-diagonal matrix:
\begin{equation}\label{eq:tri_diag}
T_k=\left(
\begin{array}{cccc}
\alpha_1 & \beta_1 & \cdots & 0 \\
\beta_1 & \alpha_2 & \ddots & \vdots\\
\vdots & \ddots & \ddots & \beta_{k-1}\\
0&\cdots & \beta_{k-1} & \alpha_k
\end{array}
\right).
\end{equation}
The Lanczos procedure is actually derived by comparing the left
and right hand sides of $WQ_k \approx Q_kT_k$ (cf. Section
\ref{sec:BLWS}).

\begin{algorithm}[h]
   \caption{The Lanczos Procedure}
   \label{Lanczos}
\begin{algorithmic}
   \STATE  {\bfseries Input:} $m\times m$ symmetric matrix $W$, $k$.
   \STATE 1. Initialization: $r_0=q_1$; $\beta_0=1$; $q_0=0$; $l=0$.
   \STATE 2. \FOR {$l=1:k-1$}
   \STATE $q_{l+1}=r_l/\beta_l$; $l=l+1$; $\alpha_l=q_l^TWq_l$;
   \STATE $r_l=Wq_l-\alpha_l q_l - \beta_{l-1}q_{l-1}$;
   \STATE $\beta_l=\|r_l\|_2$;
   \ENDFOR
   \STATE {\bfseries Output:} $Q_k=(q_1,\cdots,q_k)$ and $T_k$ as (\ref{eq:tri_diag}).
\end{algorithmic}\vspace{-0.25em}
\end{algorithm}

After the Lanczos procedure is iterated for $k-1$ times, the EVD
of $T_k$ is computed: $T_k=V_k \Lambda_k V_k^T$. Then $W\approx
(Q_kV_k)\Lambda_k (Q_kV_k)^T$. Suppose the eigenvalues in
$\Lambda_k$ is ordered from large to small. Then the $r$ largest
eigenvalues of $W$ can be approximated by the first $r$
eigenvalues in $\Lambda_k$ (called the Ritz values of $W$) and the
leading $r$ eigenvectors of $W$ can be approximated by the first
$r$ columns of $Q_kV_k$ (called the Ritz vectors of $W$).

When computing the partial SVD of a given matrix $W$, a critical
relationship between the SVD of $W$ and the EVD of the following
augmented matrix
\begin{equation}\label{eq:Expand}
\tilde{W}=\left(
\begin{array}{cc}
0 & W\\
W^T & 0
\end{array}
\right)
\end{equation}
is used. It is depicted by the following theorem \cite{Golub96}.
\begin{theorem}\label{thm:SVD_EVD}
If the SVD of an $m\times n$ $(m\leq n)$ matrix $W$ is $W=U\Sigma
V^T$, then the EVD of $\tilde{W}$ is
\begin{equation}
\tilde{W} = \tilde{U} \left(
\begin{array}{ccc}
\Sigma & 0 & 0\\
0 & -\Sigma & 0\\
0 & 0 & 0
\end{array}
\right) \tilde{U}^T,
\end{equation}
where
\begin{equation}
\tilde{U}=\dfrac{1}{\sqrt{2}}\left(
\begin{array}{ccc}
U_1 & U_1 & \sqrt{2}U_2\\
V & -V & 0
\end{array}
\right)\quad \mbox{and}\quad (U_1,U_2)=U.
\end{equation}
\end{theorem}
So by computing the EVD of $\tilde{W}$, the SVD of $W$ can be
obtained.

When computing the SVD of $W$, the Lanczos method is actually
implicitly applied to $\tilde{W}$ with the initial vector
$\tilde{q}_1$ being chosen as
\begin{equation}\label{eqn:initial_q1}
\tilde{q}_1 =(u_1^T,0^T)^T,
\end{equation}
in order to exploit the special structure of $\tilde{W}$.
Accordingly, $W$ can be approximated as $W\approx U_k B_k V_k^T$,
where columns of $U_k$ and $V_k$ are orthonormal and $B_k$ is
\emph{bi-diagonal}. Then the approximate singular values/vectors
of $W$ can be obtained after computing the SVD of $B_k$. For more
details, please refer to \cite{Propack}.

The Lanczos method has some important properties \cite{Golub96}.
First, the Ritz values of $W$ converge to the largest
eigen/singular values of $W$ quickly when $k$ grows, so do the
Ritz vectors. Second, as it only requires solving the EVD/SVD of a
relatively small and banded matrix $T_k$/$B_k$, the partial EVD/SVD is
usually faster than the full EVD/SVD when the required number $r$ of
eigen/singular vectors is relatively small (e.g., when $r <
0.25p$). Third, the Lanczos procedure terminates when an invariant
subspace is found. Fourth, the orthogonality among the columns of
$Q_k$ is easily lost when the Lanczos procedure goes on. Hence,
reorthogonalization is usually necessary when $k$ is relatively
large. Unfortunately, reorthogonalization is expensive. So PROPACK monitors
the orthogonality among the columns of $Q_k$ and only reorthogonalizes
part of the columns whose orthogonalities with other columns deteriorate.

\section{Ideas to Improve}
\label{sec:ideas} We notice that if we solve the partial SVD in
each iteration independently, the Lanczos procedure has to start
from a random initial vector $q_1$ as no apriori information is
available. Random initialization makes the size $k$ of $B_k$
unpredictable. If $q_1$ is not good, $k$ can be relatively large
in order for the Ritz values/vectors to achieve a prescribed
precision, making the partial SVD inefficient. Actually, during
the iterations of optimization, as the matrices $W_i$ and
$W_{i-1}$ are close to each other, any of the leading Ritz vectors
of $W_{i-1}$ should be good for initializing the Lanczos procedure
of $W_{i}$. However, a risk of this strategy is that the Lanczos
procedure may terminate quickly by outputting a small invariant
subspace containing the previous Ritz vector because the previous
Ritz vector is close to be a singular vector of $W_{i}$. So the
Lanczos procedure for $W_i$ may fail and has to restart with
another initial vector\footnote{Although in reality the Lanczos
procedure seldom terminates due to numerical error, our numerical
tests show that such choice of initial $q_1$ does not help
speeding up.}. Moreover, initializing with a \emph{vector} $q_1$
neglects the fact that we are actually seeking optimal singular
\emph{subspaces}, not a number of individual singular vectors. A
vector definitely cannot record all the information from the
previous principal singular subspaces (left and right). So, ideally we should use
the whole previous principal singular subspaces. This motivates us
to adopt the block Lanczos method for partial SVD, where the block
Lanczos procedure starts with the previous principal singular
subspaces.

\section{Block Lanczos with Warm Start}\label{sec:BLWS}
Again, we start with the block Lanczos with warm start (BLWS) for
partial EVD. The block Lanczos method is a natural generalization
of the standard Lanczos method by replacing the Krylov space
$K(W,q_1,k)$ with
\begin{equation}\label{eq:Lanczos-block}
\tilde{K}(W,X_1,k)=\span\{X_1,WX_1,\cdots,W^{k-1}X_1\},
\end{equation}
where $X_1$ is an orthonormal basis of an initial subspace.
Accordingly, the Lanczos procedure is generalized to the block
Lanczos procedure, which is to find an approximation of $W$:
$W\approx Q_k T_k Q_k^T$, where $T_k$ is a block tri-diagonal
matrix \cite{Golub96}:
\begin{equation}\label{eq:tri__block_diag}
T_k=\left(
\begin{array}{cccc}
M_1 & B_1^T & \cdots & 0\\
B_1 & M_2 &\ddots & \vdots\\
\vdots & \ddots & \ddots & B_{k-1}^T \\
0 &\cdots & B_{k-1} & M_k
\end{array}
\right),
\end{equation}
$Q_k=(X_1,\cdots,X_k)$, and columns of $Q_k$ are orthonormal. By
comparing the left and right hand sides of $WQ_k \approx Q_k T_k$, we
have
\begin{equation}
WX_l = X_{l-1}B_{l-1}^T + X_lM_l + X_{l+1}B_l,\quad
l=1,\cdots,k-1,
\end{equation}
where $B_0$ is defined to be 0. From the orthogonality among the
columns of $Q_k$, we have that
\begin{equation}
M_l=X_l^TWX_l,\quad l=1,\cdots,k.
\end{equation}
Moreover, if we define $R_l=WX_l - X_lM_l - X_{l-1}B_{l-1}^T$,
then $X_{l+1}B_l$ is the QR decomposition of $R_l$. This leads to
the block Lanczos procedure in Algorithm~\ref{Block_Lanczos}.

\begin{algorithm}[h]
   \caption{Block Lanczos Procedure}
   \label{Block_Lanczos}
\begin{algorithmic}
   \STATE  {\bfseries Input:} $m\times m$ symmetric matrix $W$, $m\times r$ orthogonal matrix $X_1$, $k$.
   \STATE 1. Initialization: $M_1=X_1^TWX_1$; $B_0=0$.
   \STATE 2. \FOR {$l=1:k-1$}
     \STATE $R_l=WX_l-X_lM_l-X_{l-1}B_{l-1}^T$;
     \STATE $(X_{l+1},B_l)=qr(R_l)$; (The QR decomposition)
     \STATE $M_{l+1}=X_{l+1}^TWX_{l+1}$;
   \ENDFOR
   \STATE {\bfseries Output:} $Q_k=(X_1,\cdots,X_k)$ and $T_k$ as
   in (\ref{eq:tri__block_diag}).
\end{algorithmic}
\end{algorithm}
After the approximation $W\approx Q_k T_k Q_k^T$ is obtained, one
may further compute the EVD of $T_k$: $T_k=U_k\Lambda_k U_k^T$,
where the eigenvalues $\lambda_i$ are ordered from large to small.
Then the leading $r$ eigenvalues and eigenvectors of $W$ is
approximated by $\lambda_1,\cdots,\lambda_r$, and $Q_kU_k(:,1:r)$,
respectively. The whole process is summarized in Algorithm
\ref{Block_Lanczos_EVD}.

\begin{algorithm}[h]
   \caption{Block Lanczos for Partial EVD}
   \label{Block_Lanczos_EVD}
\begin{algorithmic}
   \STATE  {\bfseries Input:} $m\times m$ symmetric matrix $W$, $m\times r$ orthogonal matrix $X_1$, $k$.
   \STATE 1. Compute $Q_k$ and $T_k$ by
   Algorithm~\ref{Block_Lanczos}.
   \STATE 2. Compute the EVD of $T_k$: $T_k=V_k\Lambda_k V_k^T$, where the eigenvalues on the diagonal of $\Lambda_k$
    are in a decreasing order.
   \STATE {\bfseries Output:} $U=Q_kV_k(:,1:r)$, $\Sigma=\Lambda_k(1:r,1:r)$.
\end{algorithmic}
\end{algorithm}

If we denote the block Lanczos for partial EVD
(Algorithm~\ref{Block_Lanczos_EVD}) as $BL\_EVD(W,X_1,k)$, then
our BLWS can be written as:
$$\mbox{(\bf BLWS)}\quad (U_i,\Sigma_i)=BL\_EVD(W_i,U_{i-1},k_i),$$
namely the principal eigen-subspace $U_{i-1}$ of the previous
iteration is used to initialize the block Lanczos procedure.

When using the block Lanczos method to compute the partial SVD of
a matrix $W$, similarly the block Lanczos procedure is applied to
$\tilde{W}$ shown in (\ref{eq:Expand}). Note that $\tilde{W}$ is
of special structure. So the block Lanczos procedure can be done
efficiently by skipping the zero sub-matrices of $\tilde{W}$. The
details are trivial. So we omit them.

With BLWS, compared with the standard Lanczos method, the risk of
terminating with a small invariant subspace is gone, and the principal
eigen-subspace can be updated more efficiently. As a result, the whole
optimization process can be sped up a lot.

\subsection{More Tricks for Acceleration}
Recall that in the standard Lanczos procedure, the orthogonality
among the columns of $Q_k$ is easily lost when $k$ grows. So is
the block Lanczos procedure. As reorthogonalization is expensive,
we further require that the number $k$ of performing the block
Lanczos procedure is small, such that reorthogonalization can be
waived. In our experiments, we typically set $k=2$, namely the
block Lanczos procedure is performed only once. Although such a
fixed and small value of $k$ cannot result in high precision
principal singular subspaces when the block Lanczos procedure is
randomly initialized, it does produce high precision principal
singular subspaces when the block Lanczos procedure is initialized
with the previous principal singular subspaces. This is because $W_i$
is close to $W_{i-1}$. So the previous principal singular subspaces
is already close to the principal singular subspaces of $W_i$. Then
the block Lanczos procedure improves them and produce better
estimated principal singular subspaces. Note that keeping $k$ small
has multiple advantages. First, it waives the necessity of
expensive reorthogonalization. Second, it saves the computation in
performing the block Lanczos procedure. Third, the SVD of $B_k$
also becomes cheap because the size of $B_k$ is small.

In the standard block Lanczos method for partial SVD, the initial
subspace is chosen as $\tilde{X}_1=(U_{i-1}^T,0)^T$ or
$\tilde{X}_1=(0,V_{i-1}^T)^T$ (cf. (\ref{eqn:initial_q1})), where
$U_{i-1}$ and $V_{i-1}$ are the estimated left and right principal
singular subspaces obtained in the previous iteration, respectively. However,
such an initialization only utilizes half of the information
provided by the previous iteration. So our BLWS technique
uses $\tilde{X}_1=\frac{1}{\sqrt{2}}(U_{i-1}^T,V_{i-1}^T)^T$ as the
initial subspace. In this way, the precision of obtained principal
singular subspaces is higher when the block Lanczos procedure is
performed for the same number of steps.

\subsection{Handling Variant Dimensions of Principal Singular Subspaces}
The above exposition assumes that the dimension $r$ of the principal singular
subspaces is known and fixed along iteration. In reality, $r$ is unknown and has
to be dynamically predicted before the partial SVD is computed \cite{Lin09,Toh09,Cai08,Ma09}.
Hence $r$ actually varies along iteration. In this case, BLWS simply
outputs Ritz values/vectors according to the predicted $r$ in the current iteration
and the block Lanczos procedure is still initialized with the
principal singular subspaces output by last iteration. We have observed that for
many NNM problems, the predicted $r$ quickly stabilizes. So variant dimensions
of principal singular subspaces at the early iterations do not affect the effectiveness of BLWS.

\section{Experimental Results}\label{sec:experiments}
Our BLWS technique is a general acceleration method. Given an
algorithm to solve a NNM problem, a user
only has to replace the SVD computation in the algorithm with BLWS
and may obtain noticeable speedup.

As examples, in this section we apply our BLWS technique to two
popular problems: Robust PCA (RPCA) \cite{Wright09} and Matrix
Completion (MC) problems \cite{Candes08}. For each problem, we
compare the original chosen algorithm and its BLWS improved
counterpart in the aspect of computation time. The accuracies of
obtained solutions are also shown in order to ensure that the
correct solutions are approached. All experiments are run on the
same workstation with two quad-core 2.53GHz Intel Xeon E5540
CPUs, running Windows Server 2008 and Matlab (Version 7.7).

For the RPCA problem, we generate the synthetic data in the same
way as that in \cite{Lin09}. Namely, $A$ is generated according to
the independent random orthogonal model \cite{Wright09}, $E$ is a
sparse matrix whose support is independent and the entry values
are uniformly distributed in $[-500,500]$, and $D=A+E$. For
simplicity, we only focus on $m\times m$ square matrices and the
parameter $\lambda$ is fixed at $1/\sqrt{m}$, as suggested by
Wright et al. \cite{Wright09}. The value of $m$ is chosen in
$\{500,1000,2000,3000\}$. The rank of $A$ is chosen as $10\%m$,
and the number of corrupted entries (i.e., $\|E\|_{l_0}$) is
$10\%m^2$. We choose the ADM method \cite{Yuan09,Lin09} to solve
the PRCA problem.

The data for the MC problem is generated in the same way as that
in \cite{Cai08}. Namely, an $m\times m$ matrix $A$ with rank $r$
is generated by first sampling two $m\times r$ factor matrices
$M_L$ and $M_R$ independently, each having i.i.d. Gaussian
entries, and then multiplying them: $A=M_LM_R^{T}$. Finally, the
set of observed entries is sampled uniformly at random. We choose
the SVT algorithm \cite{Cai08} to solve the MC problem.

Table 1 shows detailed comparison between ADM and
BLWS accelerated ADM for solving the RPCA problem. We can see that
BLWS-ADM roughly costs less than 1/3 time of ADM, achieving the
same accuracy, and the total number of iterations does not change or
only increases slightly. Similar phenomenon can also be observed in Table
2, which lists the comparison results for solving the
MC problem.

\begin{table*}
\begin{center}
\label{BLWS-ADM} \caption{BLWS-ADM vs. ADM on different synthetic
data. $\hat{A}$ and $\hat{E}$ are the computed low rank and sparse
matrices and $A$ is the ground truth. }
\begin{tabular}{|l|l|lllll|}
\hline
{$m$}    & $method$ & {$\frac{\| \hat A - A \|_F}{\|A\|_F}$} & $rank(\hat A)$ & $\|\hat E\|_{l_0}$  & \#$iter$  & $time(s)$\\
\hline
500 & ADM & 5.27e-006 & 50 & 25009 & 30 & 4.07 \\
500 & BLWS-ADM & 9.64e-006 & 50 & 25008 & 30 & 2.07 \\
\hline
1000 & ADM & 3.99e-006 & 100 & 100021 & 29 & 23.09 \\
1000 & BLWS-ADM & 6.05e-006 & 100 & 100015 & 30 & 9.25 \\
\hline
2000 & ADM & 2.80e-006 & 200 & 400064 & 28 & 154.80 \\
2000 & BLWS-ADM & 4.30e-006 & 200 & 400008 & 30 & 53.37 \\
\hline
3000 & ADM & 2.52e-006 & 300 & 900075 & 28 & 477.13 \\
3000 & BLWS-ADM & 3.90e-006 & 300 & 900006 & 30 & 149.19\\
\hline
\end{tabular}
\end{center}
\end{table*}

\begin{table*}
\begin{center}
\label{BLWS-SVT} \caption{BLWS-SVT vs. SVT on different synthetic
data. $\hat{A}$ is the recovered low rank matrix and $A$ is the
ground truth. $m$ is the size of matrix and $s$ is the number of
sampled entries. $d_r=r(2m-r)$ is the number of degrees of freedom
in an $m\times m$ matrix of rank $r$.}
\begin{tabular}{|llll|l|lll|}
\hline
{$m$} & $r$ & $s/d_r$ & $s/m^2$ & $algorithm$ & $time(s)$ & \#$iter$ & {$\frac{\| \hat A - A \|_F}{\|A\|_F}$}  \\
\hline \hline

\comment{
1000 & 10 & 6 & 0.12 & SVT & 11.81 & 117 & 1.71e-004 \\
1000 & 10 & 6 & 0.12 & BLWS-SVT & 3.34 & 118 & 1.65e-004 \\
\hline
1000 & 50 & 4 & 0.39 & SVT & 75.64 & 113 & 1.63e-004 \\
1000 & 50 & 4 & 0.39 & BLWS-SVT & 27.21 & 113 & 1.64e-004 \\
\hline

1000 & 100 & 3 & 0.57 & SVT & 160.55 & 129 & 1.69e-004\\
1000 & 100 & 3 & 0.57 & BLWS-SVT & 77.67 & 129 & 1.69e-004\\

\hline }
5000 & 10 & 6 & 0.024 & SVT & 72.57 & 123 & 1.73e-004 \\
5000 & 10 & 6 & 0.024 & BLWS-SVT & 20.02 & 123 & 1.74e-004 \\
\hline
5000 & 50 & 5 & 0.1 & SVT & 438.81 & 107 & 1.63e-004 \\
5000 & 50 & 5 & 0.1 & BLWS-SVT & 279.08 & 108 & 1.55e-004 \\
\hline

5000 & 100 & 4 & 0.158 & SVT & 1783.26 & 122 & 1.73e-004 \\
5000 & 100 & 4 & 0.158 & BLWS-SVT & 1175.91 & 122 & 1.74e-004\\

\hline
10000 & 10 & 6 & 0.012 & SVT & 135.90 & 123 & 1.68e-004 \\
10000 & 10 & 6 & 0.012 & BLWS-SVT & 42.80 & 123 & 1.70e-004 \\
\hline

10000 & 50 & 5 & 0.050 & SVT & 1156.25 & 110 & 1.58e-004\\
10000 & 50 & 5 & 0.050 & BLWS-SVT& 736.01 & 110 & 1.60e-004\\

\hline
20000 & 10 & 6 & 0.006 & SVT & 251.13 & 123 & 1.74e-004 \\
20000 & 10 & 6 & 0.006 & BLWS-SVT & 101.47 & 124 & 1.68e-004 \\
\hline
30000 & 10 & 6 & 0.004 & SVT & 449.34 & 124 & 1.75e-004 \\
30000 & 10 & 6 & 0.004 & BLWS-SVT  & 171.40 & 125 & 1.69e-004 \\
\hline

\end{tabular}
\end{center}
\end{table*}

\section{Discussions}\label{sec:discussion}
Although we have presented numerical results to testify to the
effectiveness of BLWS, currently we have not rigorously proved the
correctness of BLWS. We guess that BLWS can work well for most
NNM problems. This is due to Theorem 9.2.2
of \cite{Golub96}, which implies that when there is sufficient gap
between the $r$-th and the $(r+1)$-th eigenvalues, the errors in
the Ritz values can be well controlled. As NNM problems typically involve singular value
thresholding (\ref{eqn:SVT}), such a gap should exist. However, a
rigorous proof is still under exploration.

\section{Conclusions}\label{sec:conclusion}
In this paper, we introduce the block Lanczos with warm start
technique to accelerate the partial SVD computation in NNM problems. Both the block Lanczos procedure and the
initialization with the previous principal singular subspaces can take
full advantage of the information from last iteration. Our BLWS
technique can work in different algorithms that solve rank
minimization problems. The experimental results indicate that our
BLWS technique usually makes its host algorithm at least two to
three times faster.

%\begin{acknowledgements}
%\end{acknowledgements}

\bibliographystyle{spmpsci}      % mathematics and physical sciences
\bibliography{BL1}   % name your BibTeX data base

%\begin{thebibliography}{}
% etc
%\end{thebibliography}

\end{document}